# BIFURCATION ANALYSIS OF THE CONVERSION DEGREES IN SYSTEMS BASED ON THE CASCADE OF TANK REACTORS

*Marek Berezowski*
*Silesian University of Technology*
*Institute of Mathematics*
*44-100 Gliwice, ul. Kaszubska 23, Poland, e-mail: Marek.Berezowski@polsl.pl*

**Abstract**
The scope of the paper is the analysis of the posibility of increasing conversion in different types of systems based on continuous stirred tank reactors. A numerical model of the cassade with constant direction of material flow is tested, a model of changes of direction flow and a relaxation model. The analysis is conducted by means of parametric continuation method and numerical simulation, with the designation of fitting diagrams and time series. An essential impact of the relaxation method on the increase of the conversion degree is indicated.

**Keywords**: chemical reactors; dynamic simulation; mathematical modelling; nonlinear dynamics; numerical analysis; optimisation

## 1. Introduction

Industrial processes that involve chemical reactions often use continuous stirred tank reactors (CSTR). Such systems have been largely discussed in research and professional publications. The basic task for designers of reactor systems is to select apparatuses and their configuration to provide the highest possible conversion degree of inflowing raw material. There are many elaborations and analyses of tubular reactors with reverse flow (Sheintuch, 2005; Sheintuch & Nekhamkina, 2004; Glöckler, Kolios & Eigenberger, 2003; Jeong & Luss, 2003; Annalad, Scholts, Kuipers, & Swaaij, Part I and Part II, 2002; Khinast, Jeong & Luss, 1999; Řeháček, Kubiček & Marek, 1998; Purwono, Budman, Hudgins, Silveston & Matros, 1994; Gupta & Bhata, 1991). The analysis of the dynamics of reverse flow system used for cascade tank reactors was discussed in (Żukowski & Berezowski, 2000; Kulik & Berezowski, 2008; Berezowski & Kulik, 2009), where a possibility of the ocurrence of chaotic oscillations was proved.

According to industrial and laboratory practice, as well as to theoretical analyses, for certain parameter ranges, the cascade enables good conversion of the raw material, whereas, for other ranges very low conversion degrees are achieved, meaning that practically chemical reactions do not occur. One of the ways of increasing the conversion degree in such systems is to use cyclical reverse-flow of the raw material. Such method makes it possible to store the generated heat energy. For tubular gas reactors with fixed-bed this energy is accumulated, first and foremost, in the catalyst, whereas, in the case of liquid tank reactors the heat is collected in the reaction mix.

Likewise, it is evident that in the reverse flow system, high conversion degree may only be obtained for specific values of the system parameters, including the time between successive changes of the flow direction. The method presented in this paper enables an increment of the conversion degree also when the reverse flow renders low conversion. It involves - after the change in the direction of the feed flux - cutting off, for a certain time, the flow to the reactors, consisting of a low-concentration product mix and the simultaneous supply and receipt of the product from the apparatuses containing low concentrations of the raw material. Such time is referred to as relaxation time. Thanks to this procedure the raw material in the cut-off reactors has better chance of higher conversion. After the passage of the relaxation time, the previously cut-off reactors are again fed with fluxes and the substance subjected to conversion is collected.

The theoretical analysis of the above-mentioned systems was conducted on the bases of applicable bifurcation diagrams and time series. They concern both steady states, as well as enforced oscillation states, which are a natural consequence of reverse-flow applications. The method of designating bifurcation diagrams was discussed in (Berezowski, 2010).

## 2. Mathematical models, bifurcation diagrams and analysis of the results

The analyzed cascade consists of two identical adiabatic continuous stirred tank reactors CSTR 1 and CSTR 2 (Fig.1).

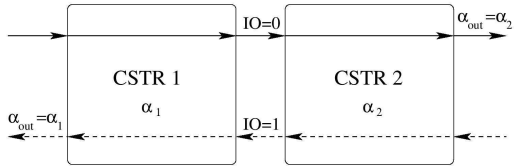

Fig. 1. Schematic diagram of the cascade of continuous stirred tank reactors with constant and reverse flow direction.

The balance equations for particular apparatuses are expressed by the following differential system:

$$\frac{d\alpha_1}{d\tau} + \alpha_1 = IO\alpha_2 + \phi(\alpha_1) \qquad (1)$$

$$\frac{d\alpha_2}{d\tau} + \alpha_2 = (1 - IO)\alpha_1 + \phi(\alpha_2) \qquad (2)$$

where $IO$ is a variable determining the direction of the reacting flux flow. Accordingly, if the flow passes from CSTR 1 to CSTR 2, $IO = 0$. Otherwise, $IO = 1$. The use of this variable will be very important in the reverse-flow system applied to the discussed cascade. It may be inferred from the adiabatic nature of the apparatusses that the temperature of the reactive mass is proportional to its concentration, which, in the assumed dimensionless notations means that $\Theta = \alpha$. Hence, it is unnecessary to introduce additional equations referring to the reacting heat into the model.

### 2.1. Model and bifurcation diagram of the cascade with constant direction of reacting flux flow

Whitin the scope of the paper, the bifurcation diagram of the cascade was determied at first, assuming that the feed flux does not change its direction. It may be inferred from the adiabatic nature of the reactors that the occurrence of autonomuos temperature and concentration oscillations is impossible. This means that the only attractors may be the steady states. On the grounds of the parametric continuation method described in (Berezowski, 2010), a diagram of the states was designated and marked with the "$ss$" line (see Figure 2).

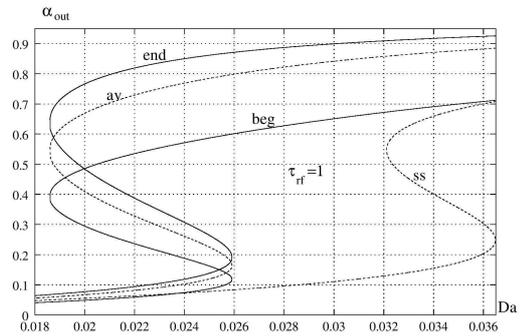

Fig. 2. Conversion degree at the cascade outlet. $\tau_{rf} = 1$.

Variable $\alpha_{out}$ is the conversion degree at the system outlet. In all caculations within the framework of the paper the following parameter vaules were used: $\gamma = 15$, $\beta = 0.65$, $n = 1.5$.

In view of the symmetry of the cascade, the calculations may be based on the assumption of any value of the variable controlling the flow direction, i.e. $IO = 0$ or $IO = 1$. In the first case, $\alpha_{out} = \alpha_2$, whereas, in the second one, $\alpha_{out} = \alpha_1$.

By equations (1)-(2), it may be concluded that in the steady state the following dependencies apply:

$$f_1 = IO\alpha_2 + \phi(\alpha_1) - \alpha_1 = 0 \quad (3)$$
$$f_2 = (1-IO)\alpha_1 + \phi(\alpha_2) - \alpha_2 = 0 \quad (4).$$

On the grounds of the parametric continuation method, the following relation is derived:

$$\overline{\alpha}_{k+1} = \overline{\alpha}_k - \overline{\overline{J}}_k^{-1}\overline{w}_k sign(\det \overline{\overline{J}})\Delta p \quad (5)$$
$$p_{k+1} = p_k + sign(\det \overline{\overline{J}})\Delta p \quad (6)$$

where $k$ is the number of the continuation code, $sign$ is the sign of martix $\overline{\overline{J}}$ determinant, $\Delta p$ an increment of the value of the continuation parameter, whereas:

$$\overline{\alpha} = \begin{bmatrix} \alpha_1 \\ \alpha_2 \end{bmatrix} \quad (7)$$

$$\overline{\overline{J}} = \begin{bmatrix} \dfrac{\partial f_1}{\partial \alpha_1} & \dfrac{\partial f_1}{\partial \alpha_2} \\ \dfrac{\partial f_2}{\partial \alpha_1} & \dfrac{\partial f_2}{\partial \alpha_2} \end{bmatrix} \quad (8)$$

$$\overline{w} = \begin{bmatrix} \dfrac{\partial f_1}{\partial p} \\ \dfrac{\partial f_2}{\partial p} \end{bmatrix} \quad (9).$$

At the same time:

$$\frac{\partial f_1}{\partial \alpha_1} = \frac{\partial \phi}{\partial \alpha_1} - 1 \quad (10)$$

$$\frac{\partial f_1}{\partial \alpha_2} = IO \quad (11)$$

$$\frac{\partial f_2}{\partial \alpha_1} = 1 - IO \quad (12)$$

$$\frac{\partial f_2}{\partial \alpha_2} = \frac{\partial \phi}{\partial \alpha_2} - 1 \quad (13).$$

Assuming that there is a single $A \to B$ reaction of the *n-th* order occuring in the reactors, the kinetic functions have the following form:

$$\phi(\alpha_i) = Da(1-\alpha_i)^n \exp\left(\gamma \frac{\beta\alpha_i}{1+\beta\alpha_i}\right); i = 1,2 \quad (14).$$

This means that the applicable derivatives are described by dependencies:

$$\frac{\partial \phi}{\partial \alpha_i} = \left[Da(1-\alpha_i)^{n-1}\left(-n + (1-\alpha_i)\frac{\gamma\beta}{(1+\beta\alpha_i)^2}\right)\right]\exp\left(\gamma\frac{\beta\alpha_i}{1+\beta\alpha_i}\right)$$
$$i = 1,2 \quad (15).$$

Assuming, for example, that bifurcation parameter $p$ is Damköhler number $Da$, vector $\overline{w}$ components are, respectively:

$$\frac{\partial f_i}{\partial p} = \frac{\partial f_i}{\partial Da} = (1-\alpha_i)^n \exp\left(\gamma\frac{\beta\alpha}{1+\beta\alpha}\right)$$
$$i = 1,2 \quad (16).$$

As indicated in Fig. 2. high conversion degrees are achieved for $Da > 0.032$, but, in this specific range the phenomenon of the multiplicity of steady states occurs.

**2.2. Model and bifurcation diagram of the cascade with reverse reacting flow**

In this section of the paper a theoretical analysis of the cascade with reverse flow is conducted. The dynamics of this process was already discussed in (Żukowski & Berezowski, 2000; Kulik & Berezowski, 2008; Berezowski & Kulik, 2009), where a possibility of the occurrence of chaotic oscillations was proved.

Assuming that the time period betweeen successive changes of the flow direction is $\tau_{rf}$, the basic oscilation period of variable $\alpha_{out}$ at the system outlet is also equal to $\tau_{rf}$. The output is CSTR 2 ($IO = 0$) and CSTR 1 ($IO = 1$), alternately.

In the course of the analysis of the system, the next two bifurcation diagrams were designated, referring, this time, to enforced cycles with time $\tau_{rf}$ (Fig. 2). To achieve this, the parametric continuation method described in (Berezowski, 2010) was reapplied. The two curves marked in Fig.2. and labelled: "*beg*" and "*end*" refer to the value of $\alpha_{out}$ at the beginning and at the end of the cycle. Curve "*av*" is the average value of $\alpha_{out}$, calculated in a given cycle.

After considering periodic change of the raw material flux flow direction in model (1)-(2), the general balance equations have the following discreted form:

$$\frac{d\alpha_{1,j+1}}{d\tau} + \alpha_{1,j+1} = IO\alpha_{2,j+1} + \phi(\alpha_{1,j+1}) \quad (17)$$

$$\frac{d\alpha_{2,j+1}}{d\tau} + \alpha_{2,j+1} = (1-IO)\alpha_{1,j+1} + \phi(\alpha_{2,j+1}) \quad (18)$$

with boundary conditions:

$$\alpha_{1,j+1}(\tau_{rf}) = \psi(\alpha_{1,j+1}(0), IO\alpha_{2,j+1}(0)) = \psi(\alpha_{1,j}(\tau_{rf}), IO\alpha_{2,j}(\tau_{rf})) \quad (19)$$

$$\alpha_{2,j+1}(\tau_{rf}) = \psi(\alpha_{2,j+1}(0), (1-IO)\alpha_{1,j+1}(0)) = \psi(\alpha_{2,j}(\tau_{rf}), (1-IO)\alpha_{1,j}(\tau_{rf})) \quad (20)$$

where $j$ is the number of the cycle, whereas $\psi$ is an integral transformation of equations (17) and (18).

In accordance with the previous assumptions, variable $IO = 0$ when $j$ is an even number, but $IO = 1$ when $j$ is an odd number. The conversion degree at the outlet of the system is decsribed by the equation:

$$\alpha_{out,j} = IO\alpha_{1,j} + (1-IO)\alpha_{2,j} \quad (21).$$

If the sought solutions are fixed oscilation points enforced in time $\tau_{rf}$, the following relations hold for each cycle:

$$\alpha_1(\tau_{rf}) = \alpha_2(0) \quad (22)$$
$$\alpha_2(\tau_{rf}) = \alpha_1(0) \quad (23).$$

Hence, the parametric continuation method should be used for the system of equations:

$$F_1 = \alpha_1(\tau_{rf}) - \alpha_2(0) = 0 \quad (24)$$
$$F_2 = \alpha_2(\tau_{rf}) - \alpha_1(0) = 0 \quad (25)$$

which, subsequently, leads to the following recurrence process:

$$\overline{\alpha}_{k+1}(\tau_{rf}) = \overline{\alpha}_k(\tau_{rf}) - \overline{\overline{A}}_k^{-1}\overline{b}_k sign(\det \overline{\overline{A}})\Delta Da \quad (26)$$

$$Da_{k+1} = Da_k + sign(\det \overline{\overline{J}})\Delta Da \quad (27).$$

where $k$ is the continuation step number, *sign* is the sign of martix $\overline{\overline{J}}$ determinant, $\Delta Da$ is the increment of the value of Damköhler number, while:

$$\overline{\alpha}(\tau_{rf}) = \begin{bmatrix} \alpha_1(\tau_{rf}) \\ \alpha_2(\tau_{rf}) \end{bmatrix} \quad (28)$$

$$\overline{\overline{A}} = \begin{bmatrix} \frac{\partial F_1}{\partial \alpha_1(0)} & \frac{\partial F_1}{\partial \alpha_2(0)} \\ \frac{\partial F_2}{\partial \alpha_1(0)} & \frac{\partial F_2}{\partial \alpha_2(0)} \end{bmatrix} = \begin{bmatrix} \frac{\partial \alpha_1(\tau_{rf})}{\partial \alpha_1(0)} & -1 \\ \frac{\partial \alpha_2(\tau_{rf})}{\partial \alpha_1(0)} - 1 & \frac{\partial \alpha_2(\tau_{rf})}{\partial \alpha_2(0)} \end{bmatrix} \quad (29)$$

$$\overline{b} = \begin{bmatrix} \frac{\partial F_1}{\partial Da} \\ \frac{\partial F_2}{\partial Da} \end{bmatrix} = \begin{bmatrix} \frac{\partial \alpha_1(\tau_{rf})}{\partial Da} \\ \frac{\partial \alpha_2(\tau_{rf})}{\partial Da} \end{bmatrix} \quad (30).$$

The partial derivatives in the above equations may be easily calculated by numerical methods. The designated $\alpha_1(\tau_{rf})$ and $\alpha_2(\tau_{rf})$ are the initial and final values of conversion degree $\alpha_{out}$ in a given cycle. In Fig.2. they were marked as the lines: "*beg*" and "*end*".

In view of process considerations, the average value of the conversion degree at the outlet, calculated as the integral average in a given cycle (line "*av*" in Fig.2) is essential. In comparison with the conversion degree in the cascade with constant flow direction, it is obvious that in the case of reverse-flow, the conversion of the system was considerably increased. Range $0.0185 < Da < 0.032$ is of special importance here, because within this range the cascade without the reverse-flow renders only low conversion degrees. On the other hand, the system with the reverse flow renders significantly higher conversion in the same range. Accordingly, for *Da*=0.028 the use of the reverse flow increases the average conversion by almost 800%. It should also be emphasised that in range $0.0185 < Da < 0.026$ the solutions are characterised by the multiplicity of oscilation states.

In Fig.3. exemplary time series of conversion degree $\alpha_{out}$ for $\tau_{rf} = 1$ and *Da*=0.022 were plotted.

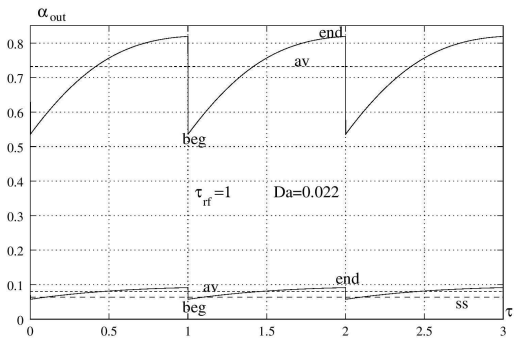

Fig. 3. Time series of the conversion degree at the cascade outlet. $\tau_{rf} = 1$.

The curves at the top and bottom of the graph refer, respectively, to the upper and lower oscilation state shown in Fig.2. Furthermore, curve "*ss*" in Fig.3. refers to the conversion degree of the cascade with constant flow direction.

Conversely, in Fig. 4 the bifurcation diagrams plotted for $\tau_{rf} = 6$ are shown.

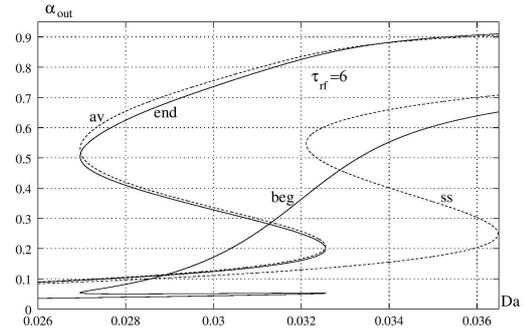

Fig. 4. Bifurcation diagrams of the conversion degree at the cascade outlet. $\tau_{rf} = 6$.

Also in this case, a significant impact of the reverse flow on the conversion degree is clearly observable. A phenomenon chracteristic for this case is that curve "*av*" is above "*end*" curve. The explanation of this phenomenon is presented by average of the time series in Fig.5, where typical extremes are evident.

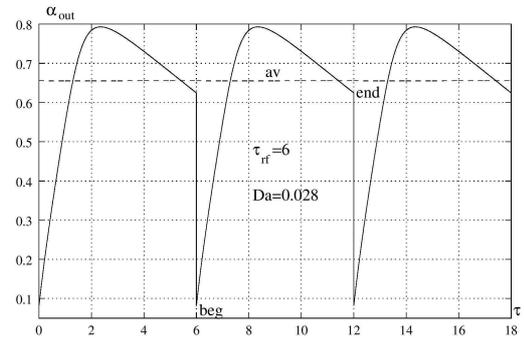

Fig. 5. Time series of the conversion degree at the cascade outlet. $\tau_{rf} = 6$.

## 2.3. Cascade with the changes of the flow direction of the reacting flux and relaxation time

As indicated in Fig.4, for $Da < 0.027$ the discussed cascade renders only low conversion degrees, despite constant or reverse material flow direction. However, it is possible to significantly increase the average conversion of the raw material. To achieve this, after changing the flow direction, the inflow of the flux to the outlet reactor should be cut off for a certain time, referred to as relaxation time $\tau_{rel} \leq \tau_{rf}$, so that the product is collected only from the reactor into which the raw material is currently fed (Fig.6 – case $IO = 1$).

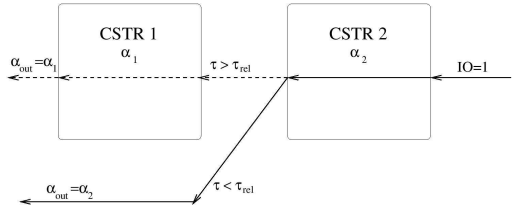

Fig. 6. Schematic diagram of the cascade with relaxation.

Thanks to such procedure, the raw material accumulated in the first reactor will be subjected to conversion, whereas the product amassed in the second reactor will be collected. The value of time $\tau_{rel}$ should be selected in a manner that secures the highest average value of variable $\alpha_{out}$.

In Fig.7 the time series for variations of $\alpha_{out}$ were plotted for $\tau_p = 6$ and $Da = 0.0265$, and the use of relaxation time $\tau_{rel} = 4.5$, which, for this particular case is the optimal value, i.e. the value that guarantees the highest average values.

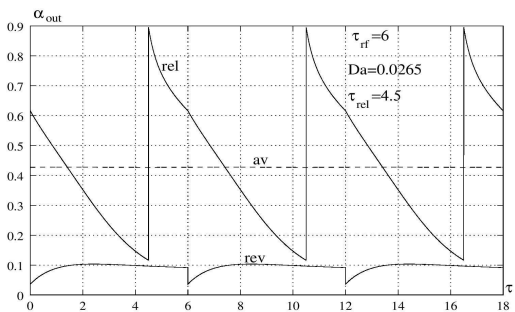

Fig. 7. Time series of the conversion degree at the relaxation cascade outlet.
$\tau_{rf} = 6$, $\tau_{rel} = 4.5$.

Curve "*rev*" refers to the cascade with the reverse flow, "*rel*" to the relaxation system, whereas line "*av*" is the average value of variable $\alpha_{out}$ after the application of the relaxation procedure. It follows from Fig.7. that this average value was increased by about 400 % in comparison with the conversion oferred by the previously discussed systems. It should also be noted that, in some cases, the optimal relaxation time may be equal to the time of cyclical changes in the raw material flow direction, namely: $\tau_{rel} = \tau_p$.

### 3. Conclusions

The theoretical analysis of systems based on the cascades of continuous stirred tank reactors included three systems: the cascade with constant raw material flow direction, the cascade with reverse flow and the cascade in which both reverse flow and relaxation of the reacting mass were used. The analysis was based on the bifurcation diagrams and time series of steady states and enforced oscilations. A positive impact of the reverse flow and relaxation on the conversion degree at the cascade outlet was proved. Likewise, a possibility of the occurrence of multiple cyclical states in the cascade with the enforced reverse flow was indicated.

### Notations

| | |
|---|---|
| $c_p$ | heat capacity, *kJ/(kg K)* |
| $C_A$ | concentration of component A, $kmol/m^3$ |
| $Da$ | Damköhler number $\left( = \dfrac{V_R(-r_0)}{\dot{F} C_{A0}} \right)$ |
| $E$ | activation energy, *kJ/kmol* |
| $\dot{F}$ | volumetric flow rate, $m^3/s$ |
| $(-\Delta H)$ | heat of reaction, *kJ/kmol* |
| $k$ | reaction rate constant, $1/s(m^3/kmol)^{n-1}$ |
| $n$ | order of reaction |
| $(-r)$ | rate of reaction, $(= kC^n)$, $kmol/(m^3 s)$ |
| $R$ | gas constant, *kJ/(kmol K)* |

| | | | |
|---|---|---|---|
| $t$ | time, $s$ | $\Theta$ | dimensionless temperature $\left(=\dfrac{T-T_0}{\beta T_0}\right)$ |
| $T$ | temperature, $K$ | | |
| $V$ | volume, $m^3$ | | |

*Greek letters*

$\alpha$    degree of conversion $\left(=\dfrac{C_{A0}-C_A}{C_{A0}}\right)$

$\beta$    dimensionless number related to adiabatic temperature increase $\left(=\dfrac{(-\Delta H)C_{A0}}{T_0 \rho c_p}\right)$

$\gamma$    dimensionless number related to activation energy $\left(=\dfrac{E}{RT_0}\right)$

$\rho$    density $\left(=\dfrac{kg}{m^3}\right)$

$\tau$    dimensionless time $\left(=\dfrac{\dot{F}}{V_R}t\right)$

$\tau_{rf}$    dimensionless reverse period

$\tau_{rel}$    dimensionless relaxation time

*Subscripts*
0    refers to feed
1,2    refers to number of reactor
R    refers to reactor

## References

Annaland, M., Scholts, H.A.R., Kuipers, J.A.M. & Swaaij, W.P.M. (2002). A novel reverse flow reactor coupling endothermic and exothermic reactions. Part I: comparison of reactor configurations for irreversible endothermic reactions. *Chemical Engineering Science,* 57, 833-854.

Annaland, M., Scholts, H.A.R., Kuipers, J.A.M. & Swaaij, W.P.M. (2002). A novel reverse flow reactor coupling endothermic and exothermic reactions. Part II: Sequential reactor configuration for reversible endothermic reactions. *Chemical Engineering Science,* 57, 855-872.

Berezowski, M. (2010). The application of the parametric continuation method for determining steady state diagrams in chemical engineering. *Chemical Engineering Science,* 65/19, 5411-5414.

Berezowski, M. & Kulik, B. (2009). Periodicity of chaotic solutions of the model of thermally coupled cascades of chemical tank reactors with flow reversal. *Chaos, Solitons & Fractals*, 40/1, 331-336.

Glöckler, B., Kolios, G. & Eigenberger, G. (2003). Analysis of a novel reverse-flow reactor concept for autothermal methane steam reforming. *Chemical Engineering Science,* 58, 593-601.

Gupta, V.K. & Bhata, S.K. (1991). Solution of cyclic profiles in catalytic reactor operation with periodic flow reversal. *Comp. Chem. Engng.*, 15/4, 229-237.

Jeong, Y.O. & Luss, D. (2003). Pollutant destruction in a reverse-flow chromatographic reactor. *Chemical Engineering Science,* 58, 1095-1102.

Knihast, J., Jeong, Y.O. & Luss, D. (1999). Dependence of cooled reverse-flow reactor dynamics on reactor model. *AIChE J.,* 45/2, 299-309.

Kulik, B. & Berezowski, M. (2008). Chaotic dynamics of coupled cascades of tank reactors with flow reversal, *Chemical and Process Engineering,* 29, 465-471.

Purwono, S., Budman, H., Hudgins, R.R., Silveston, P.L. & Matros, Yu.Sh. (1994). Runway in packed bed reactors operating with periodic flow reversal. *Chemical Engineering Science,* 49/24B, 5473-5487.

Řeháček, J., Kubiček, M. & Marek. M. (1998). Periodic, quasiperiodic and chaoic spatiotemporal patterns in a tubular catalytic reactor with periodic flow reversal. *Comp. Chem. Engng*., 22/1-2, 283-297.

Sheintuch, M. (2005). Analysis of design sensitivity of flow-reversal reactors: Simulations, approximations and oxidations experiments. *Chemical Engineering Science,* 60, 2991-2998.

Sheintuch, M. & Nekhamkina, O. (2004). Comparison of flow-reversal, internal-recirculation and loop reactors. *Chemical Engineering Science*, 59, 4065-4072.

Żukowski, W. & Berezowski, M. (2000). Generation of chaotic oscillations in a system with flow reversal. *Chemical Engineering Science*, 55, 339-343.